\documentstyle[amssymb]{amsart}

\theoremstyle{definition}
\newtheorem{obs}{Observation}
\newtheorem{prop}{Proposition}
\theoremstyle{remark}
\newtheorem{rem}{Remark}

\begin{document}

\baselineskip = 15 pt

\begin{titlepage}

\title[Free boson representation of Yangian double]
{\bf $q$-affine-Yangian double correspondence
\\and free boson representation of \\Yangian double with arbitrary level}

\vspace{1cm}
\author{Bo-Yu Hou \hspace{1cm}  Liu Zhao\\
Institute of Modern Physics, Northwest University, Xian 710069,
China\\
Xiang-Mao Ding\\
Institute of Theoretical Physics,
Academy of China, Beijing 100080, China
}

\date{January 1997}
\maketitle

\vspace{2cm}

\begin{abstract}
We derive a free boson representation of the Yangian double
$DY_\hbar(sl_N)$ with arbitrary level $k$ using the observation
that there is a correspondence between the $q$-affine algebra
and Yangian double associated with the same Cartan matrix.
Vertex operator and screening currents are not obtained in
the same way.
\end{abstract}

\end{titlepage}

\section{Introduction}
$q$-algebra and Yangian were proposed by Drinfeld as generalizations
of classical Lie algebras with nontrivial Hopf algebra structures
\cite{D1,D2,D3,Dr:new}.
Following the Faddeev-Reshetikhin-Takhtajan formalism
\cite{FRT}, both kinds of
algebras can be considered as associative algebras defined through
the Yang-Baxter relation (i.e. $RLL$-relations) with the structure constants
determined by the solutions of the quantum Yang-Baxter
equation (QYBE). $q$-affine algebra \cite{Dr:new}
and Yangian double \cite{KT,K,KL} with
center are respectively affine extensions of $q$-algebra and Yangian.
$q$-affine algebra corresponds to the trigonometric solution of
QYBE and Yangian corresponds to the rational one--if one considers
the Reshetikhin-Semenov-Tian-Shansky realization \cite{RS}
which is the affine analog
of Faddeev-Reshetikhin-Takhtajan formalism--and they both were
proved to have important applications in certain physical problems,
especially in describing the dynamical symmetries and calculating
the correlation functions and/or form factors of some two-dimensional
exactly solvable lattice statistical model and (1+1)-dimensional
completely integrable quantum field theories \cite{BL,LS,S2}
In such applications, the
infinite-dimensional representations of $q$-affine algebra and
Yangian double are often required, especially the representations
with higher ($k>1$) level.

In practice, realization of complicated algebra in terms of a
relatively simple one is proved to be quite effective and useful.
In this aspect, the Heisenberg algebra (or free boson) representation
has become a common method for obtaining representations of
($q$-)affine algebras. For examples, the free boson representations
of $U_q(\widehat{sl_2})$ with an arbitrary level have been obtained
in Refs.\cite{Srs,M1,M2,kim,abg}.
Free boson representation of $U_q(\widehat{sl_N})$
with level 1 was constructed in \cite{FJ}. Free boson representations
of $U_q(\widehat{sl_3})$ and $U_q(\widehat{sl_N})$ with arbitrary level
were constructed in \cite{sl3} and \cite{sln} respectively.
For the Yangian doubles,
the free field representation of $DY_\hbar(sl_2)$ with level $k$ was
constructed in \cite{konno}. The level 1 free boson representation of
$DY_\hbar(sl_N)$ was given in \cite{iohara}. But free field representations for
Yangian doubles of higher rank and with arbitrary level are still unknown.

In this paper we shall address the problem of free field representation
of the Yangian double $DY_\hbar(sl_N)$ with arbitrary level $k$.
For this purpose we largely rely on the result of \cite{sln} on
free field representation of $U_q(\widehat{sl_N})$ with arbitrary level
and observe that there is a simple correspondence between the
$q$-affine algebra $U_q(\widehat{sl_N})$ and the Yangian double
$DY_\hbar(sl_N)$. This correspondence makes our derivation of
free field representation for $DY_\hbar(sl_N)$ greatly simplified.
However, we have been unable to obtain the vertex operators
and screening currents for $DY_\hbar(sl_N)$ following the same spirit.

\section{$q$-affine and Yangian double correspondence}

In this section we first establish the correspondence
between the $q$-affine algebra $U_q(\widehat{sl_N})$ and the
Yangian double $DY_\hbar(sl_N)$. For this and the subsequent purposes
we use the Drinfeld current realizations for both algebras. Other
realizations such as the Reshetikhin-Semenov-Tian-Shansky realization
of them can be found in \cite{RS} and \cite{iohara} (which are
actually the quotient algebras of $U_q(\widehat{gl_N})$ and
$DY_\hbar(gl_N)$ respectively with respect to a Heisenberg subalgebra),
and the equivalence (algebra isomorphism) to Drinfeld
realizations were given in
\cite{DF} and \cite{iohara} respectively.
We remark that the Ding-Frenkel isomorphism only provides
an algebra isomorphism but not a Hopf algebra
isomorphism for $q$-affine algebras at least in the $\widehat{gl_2}$
case \cite{BSnew}. Whether this is also the case for Yangian doubles
is an interesting open problem.

\subsection{Drinfeld currents realization of $U_q(\widehat{sl_N})$}

$U_q(\widehat{sl_N})$ is an associative algebra generated by the
Drinfeld generators $E^{\pm,i}_n~(n\in {\Bbb Z})$, $H^i_n~(n\in {\Bbb Z})~
(i=1,~2,~...,~N-1)$ and the center $\gamma$. Let

\begin{displaymath}
K_i = \mbox{exp}\left((q-q^{-1}) \frac{1}{2} H^i_0 \right),
\end{displaymath}

\noindent then we can write the Drinfeld currents in the form of formal
power series of the complex parameter $z$ with coefficients given by the
above generators,

\begin{eqnarray}
& & H^i(z) = \sum_{n\in {\Bbb Z}}H^i_n z^{-n-1},~~
E^{\pm,i}(z) = \sum_{n \in {\Bbb Z}} E^{\pm,i}_n z^{-n-1},\nonumber\\
& & \psi^i_\pm(z) = \sum_{n \in {\Bbb Z}} \psi^i_{\pm,n} z^{-n}
\equiv K_i^{\pm 1} \mbox{exp} \left( \pm (q-q^{-1})
\sum_{\pm n > 0} H^i_n z^{-n} \right). \nonumber
\end{eqnarray}

The generating relations for $U_q(\widehat{sl_N})$ in terms of these
currents can be written as follows \cite{sln},

\begin{eqnarray}
& & [ \psi^i_\pm(z),~\psi^j_\pm(w) ] =0, \label{1} \\
& & (z-q^{a_{ij}} \gamma^{-1} w) (z-q^{-a_{ij}} \gamma w)
\psi^i_+(z) \psi^j_-(w) \nonumber\\
& &~~~~= (z-q^{a_{ij}} \gamma w) (z-q^{-a_{ij}} \gamma^{-1} w)
\psi^i_-(w) \psi^j_+(z), \label{2}\\
& & (z-q^{\pm a_{ij}} \gamma^{\mp \frac{1}{2}} w)
\psi^i_+(z) E^{\pm,j} (w) =
(q^{\pm a_{ij}} z- \gamma^{\mp \frac{1}{2}} w)
E^{\pm,j} (w) \psi^i_+(z), \label{3} \\
& & (z-q^{ \pm a_{ij}} \gamma^{\mp \frac{1}{2}} w)
E^{\pm,j} (z) \psi^i_-(w) =
(q^{ \pm a_{ij}} z- \gamma^{\mp \frac{1}{2}} w)
\psi^i_-(w) E^{\pm,j} (z), \label{4} \\
& & [ E^{+,i}(z),~E^{-,j}(w) ] = \frac{\delta^{ij}}{(q-q^{-1})zw}
\left( \delta(z^{-1}w \gamma) \psi^i_+( \gamma^{\frac{1}{2}}w)
- \delta(z^{-1}w \gamma^{-1}) \psi^i_-( \gamma^{- \frac{1}{2}}w)
\right), \label{5} \\
& & (z- q^{\pm a_{ij}} w ) E^{\pm,i}(z)E^{\pm,j}(w)
= (q^{\pm a_{ij}} z - w ) E^{\pm,j}(w)E^{\pm,i}(z), \label{6} \\
& & E^{\pm,i}(z)E^{\pm,j}(w) = E^{\pm,j}(w)E^{\pm,i}(z)
~~\mbox{for}~a_{ij}=0, \label{7} \\
& & E^{\pm,i}(z_1) E^{\pm,i}(z_2) E^{\pm,j}(w)
-(q+q^{-1}) E^{\pm,i}(z_1) E^{\pm,j}(w) E^{\pm,i}(z_2) \nonumber\\
& &~~~~~+ E^{\pm,j}(w) E^{\pm,i}(z_1) E^{\pm,i}(z_2) +
(\mbox{replacement:}~z_1 \leftrightarrow z_2) = 0
~\mbox{for}~a_{ij} = -1, \label{8}
\end{eqnarray}

\noindent where $a_{ij}$ are elements of the Cartan matrix of the type
$A_{N-1}$ and

\begin{eqnarray*}
\delta(x)= \sum_{n\in {\Bbb Z}} x^n.
\end{eqnarray*}

In this paper we only consider $q$-affine algebra and Yangian double as
associative algebras and do not care about the Hopf algebra aspect.

\subsection{The Yangian double $DY_\hbar(sl_N)$}

As an associative algebra, the Yangian double $DY_\hbar(sl_N)$ is
generated by the Drinfeld generators $\{h_{il},~e^{\pm}_{il} |
i=1,~2,~...,~N-1;~l \in {\Bbb Z}_{ \geq 0} \}$ and the center $c$. In terms of the
formal power series (Drinfeld currents)

\begin{eqnarray*}
& & H^{+}_i(u) = 1 + \hbar \sum_{l \geq 0} h_{il} u^{-l-1},~
H^{-}_i(u) = 1 - \hbar \sum_{l < 0} h_{il} u^{-l-1}, \\
& & E^{\pm}(u) = \sum_{l \in {\Bbb Z}} e^{\pm}_{il} u^{-l-1}
\end{eqnarray*}

\noindent we can write the generating relations for $DY_\hbar(sl_N)$
as follows \cite{iohara},

\begin{eqnarray}
& & [ H_i^\pm(u),~H_j^\pm(v) ] =0, \label{y1}\\
& & (u_\mp-v_\pm + B_{ij} \hbar) (u_\pm-v_\mp - B_{ij} \hbar)
H_i^+(u) H_j^-(v) \nonumber\\
& &~~~~= (u_\mp-v_\pm - B_{ij} \hbar) (u_\pm-v_\mp + B_{ij} \hbar)
H_j^-(v) H_i^+(u), \label{y2} \\
& & (u_\pm-v \mp B_{ij} \hbar) H_i^+(u) E^{\pm}_j(v)
= (u_\pm-v \pm B_{ij} \hbar) E^{\pm}_j(v) H_i^+(u), \label{y3} \\
& & (u_\mp-v \mp B_{ij} \hbar) H_i^-(u) E^{\pm}_j(v)
= (u_\mp-v \pm B_{ij} \hbar) E^{\pm}_j(v) H_i^-(u), \label{y4} \\
& & (u-v \mp B_{ij} \hbar) E^{\pm}_i(u) E^{\pm}_j(v)
=  (u-v \pm B_{ij} \hbar) E^{\pm}_j(v) E^{\pm}_i(u), \label{y5} \\
& & [ E^+_i(u),~E^-_j(v) ] =
\frac{1}{\hbar} \delta_{ij} \left(
\delta( u_- - v_+) H^+_i(v_+) - \delta( u_+ - v_-) H^-_i(v_-) \right),
\label{y6} \\
& & E^{\pm}_i(u_1) E^{\pm}_i(u_2) E^{\pm}_j(v)
-2 E^{\pm}_i(u_1) E^{\pm}_j(v) E^{\pm}_i(u_2) \nonumber\\
& &~~~~~+ E^{\pm}_j(v) E^{\pm}_i(u_1) E^{\pm}_i(u_2) +
(\mbox{replacement:}~u_1 \leftrightarrow u_2) = 0
~\mbox{for}~|i-j|=1, \label{y7}\\
& & E^{\pm}_i(u) E^{\pm}_j(v)
= E^{\pm}_j(v) E^{\pm}_i(u) ~~\mbox{for}~|i-j|>1, \label{y8}
\end{eqnarray}

\noindent where

\begin{displaymath}
u_\pm = u \pm \frac{1}{4} \hbar c
\end{displaymath}

\noindent and

\begin{displaymath}
B_{ij} = \frac{1}{2} a_{ij}.
\end{displaymath}

\subsection{$q$-affine-Yangian double correspondence}

Our central goal is to establish a free boson representation of the
Yangian double $DY_\hbar(sl_N)$. For this we would like to
use the known results \cite{sln} for the $q$-affine algebra
$U_q(\widehat{sl_N})$ by establishing a correspondence principle
between these two algebras. Such a correspondence principle
has been expected for some time and was ``quite mysterious'' as
stated in Ref.\cite{iohara}.

For the present authors, however, such a correspondence is rather
obvious by making use of the Drinfeld current realizations for both
$U_q(\widehat{sl_N})$ and $DY_\hbar(sl_N)$. For other realizations
no such an obvious observation could be obtained. We give the following

\begin{obs} \label{ob1}
($q$-affine-Yangian double correspondence).
The following gives a simple correspondence between $U_q(\widehat{sl_N})$
and $DY_\hbar(sl_N)$ as associative algebras

\begin{eqnarray*}
& & q \rightarrow \mbox{e}^{\frac{\hbar}{2}},~~~
\gamma \rightarrow \mbox{e}^{\frac{\hbar c}{2}},\\
& & z \rightarrow \mbox{e}^{u} ,\\
& & \psi^i_{\pm}(z) \rightarrow H^{\pm}_i(u),\\
& & z E^{\pm,i}(z) \rightarrow E^{\pm}_i( u)
\end{eqnarray*}

\noindent in the limit $\hbar \rightarrow 0,~u \rightarrow 0$ up to the linear
approximation in $\hbar$ and $u$.
\end{obs}

We remark that the above observation only gives a rule for obtaining
equations (\ref{y1}-\ref{y8}) from (\ref{1}-\ref{8}) and does not
imply any more fundamental
Hopf algebraic or algebraic relations.

\section{Free boson representation of $DY_\hbar(sl_N)$ with arbitrary level}

In this section we shall consider our central problem--the establishment
of a free boson representation of $DY_\hbar(sl_N)$ with arbitrary level.
For $N=2$ this problem has already been solved in Ref.
\cite{konno}. For generic $N$,
the desired expressions are rather complicated and our construction
depend largely on the observation \ref{ob1} and the result of \cite{sln}.
One crucial difference of our construction from the one in
\cite{sln} is that, in our case,
the Yangian double $DY_\hbar(sl_N)$ should be realized through
{\em ordinary} Heisenberg algebras (i.e. {\em without} deformation), whereas
in Ref.\cite{sln}, $U_q(\widehat{sl_N})$ was realized via
a set of $q$-deformed Heisenberg algebras. Therefore our observation
\ref{ob1} has to be used in somewhat a nontrivial way (for example,
the vertex operators and screening currents cannot be obtained
using our correspondence principles).

\subsection{Free bosons and Fock space}

We introduce the following set of $N^2-1$ Heisenberg algebras with generators
$a^i_n~(1 \leq i \leq N-1),~b^{ij}_n~\mbox{and}~c^{ij}_n~
(1 \leq i < j \leq N)$ with $ n \in {\Bbb Z} - \{ 0 \}$ and
$p_{a^i},~q_{a^i}~(1 \leq i \leq N-1),~p_{b^{ij}},~q_{b^{ij}},~
p_{c^{ij}},~q_{c^{ij}}~(1 \leq i < j \leq N)$,

\begin{eqnarray*}
\begin{array}{ll}
$$[ a^i_n,~a^j_m ] = (k+g) B_{ij} n \delta_{n+m,0},$$ &
$$[ p_{a^i},~q_{a^j} ] = (k+g) B_{ij},$$ \cr
$$[ b^{ij}_n,~b^{i'j'}_m ] = - n \delta^{i,i'} \delta^{i,j'}
\delta_{n+m,0},$$ &
$$[ p_{b^{ij}},~q_{b^{ij}} ] = -\delta^{i,i'} \delta^{i,j'},$$ \cr
$$[ c^{ij}_n,~c^{i'j'}_m ] =  n \delta^{i,i'} \delta^{i,j'}
\delta_{n+m,0},$$ &
$$[ p_{c^{ij}},~q_{c^{ij}} ] = \delta^{i,i'} \delta^{i,j'},$$
\end{array}
\end{eqnarray*}

\noindent where $g=N$ is the dual Coexter number for the Cartan matrix
of type $A_{N-1}$.

The Fock space corresponding to the above Heisenberg algebras can be
specified as follows. Let $| 0 \rangle$ be the vacuum state defined by

\begin{eqnarray*}
& & a^i_n| 0 \rangle = b^{ij}_n| 0 \rangle = c^{ij}_n| 0 \rangle
=0~ ( n >0),\\
& & p_{a^i}| 0 \rangle = p_{b^{ij}}| 0 \rangle = p_{c^{ij}}| 0 \rangle
=0.
\end{eqnarray*}

\noindent Define

\begin{eqnarray*}
& & | l_a,l_b,l_c \rangle = \\
& &~~~~ \mbox{exp} \left(
\sum_{i,j=1}^{N-1} \sum_{n>0} l_{a^i} \frac{1}{k+g} (B^{-1})^{ji} a^j_n
- \sum_{1 \leq i < j \leq N} l_{b^{ij}} q_{b^{ij}}
+ \sum_{1 \leq i < j \leq N} l_{c^{ij}} q_{c^{ij}} \right) | 0 \rangle,
\end{eqnarray*}

\noindent it can be shown that the following equations hold,

\begin{eqnarray*}
& & a^i_n| l_a,l_b,l_c \rangle = b^{ij}_n| l_a,l_b,l_c \rangle
= c^{ij}_n| l_a,l_b,l_c \rangle =0 ~( n > 0),\\
& & p_{a^i} | l_a,l_b,l_c \rangle = l_{a^i}  | l_a,l_b,l_c \rangle,\\
& & p_{b^{ij}} | l_a,l_b,l_c \rangle = l_{b^{ij}}  | l_a,l_b,l_c \rangle,\\
& & p_{c^{ij}} | l_a,l_b,l_c \rangle = l_{c^{ij}}  | l_a,l_b,l_c \rangle.
\end{eqnarray*}

\noindent The Fock space ${\cal F}(l_a,l_b,l_c)$ is then generated
by the actions of the negative modes of $a^i,~b^{ij},~c^{ij}$. We shall
see later that this Fock space actually forms a (Wakimoto-like \cite{wakimoto,fff})
module for the Yangian double $DY_\hbar(sl_N)$ with level $k$.

For $X= a^i,~b^{ij},~c^{ij}$, let us now define

\begin{eqnarray*}
& & X(u;A,B)=\sum_{n>0} \frac{X_{-n}}{n} (u+A \hbar)^n  -
\sum_{n>0} \frac{X_{n}}{n} (u+B \hbar)^{-n} + \mbox{log}(u+B \hbar) p_X
+q_X, \\
& & X_+(u;B)= - \sum_{n>0} \frac{X_{n}}{n} (u+B \hbar)^{-n}
+ \mbox{log}(u+B \hbar) p_X,\\
& & X_-(u;A)=\sum_{n>0} \frac{X_{-n}}{n} (u+A \hbar)^n
+q_X,\\
& & X(u;A) = X(u;A,A),~~X(u)=X(u,0).
\end{eqnarray*}

\noindent Then we have

\begin{displaymath}
: \mbox{exp} \left( X(u;A,B) \right) :
= \mbox{exp} \left( X_-(u;A) \right) \mbox{exp} \left( X_+(u;B) \right).
\end{displaymath}

\noindent Following the standard quantum field theory we have

\begin{equation}
X^{\alpha}(u;A,B) X^{\beta}(v;C,D) =
\langle X^{\alpha}(u;A,B) X^{\beta}(v;C,D) \rangle
+ :X^{\alpha}(u;A,B) X^{\beta}(v;C,D):, \label{contract}
\end{equation}

\noindent where
\footnote{Here and below, all OPE relations should be understood
to hold in the analytic continuation sense.}

\begin{eqnarray*}
& & \langle a^{i}(u;A,B) a^{j}(v;C,D) \rangle
=(k+g) B_{ij} \mbox{log}(u-v+(B-C) \hbar),\\
& & \langle b^{ij}(u;A,B) b^{i'j'}(v;C,D) \rangle
=- \delta^{ii'} \delta^{jj'} \mbox{log}(u-v+(B-C) \hbar),\\
& & \langle c^{ij}(u;A,B) c^{i'j'}(v;C,D) \rangle
= \delta^{ii'} \delta^{jj'} \mbox{log}(u-v+(B-C) \hbar),
\end{eqnarray*}

\noindent and all other contractions vanish. From eq.(\ref{contract}) it is
easy to calculate that

\begin{eqnarray}
& & :\mbox{exp} \left( X^{\alpha}(u;A,B) \right):
:\mbox{exp} \left( X^{\beta}(v;C,D) \right) :  \nonumber \\
& &~~~~= \mbox{exp} \left( \langle X^{\alpha}(u;A,B)
X^{\beta}(v;C,D) \rangle \right)
:\mbox{exp} \left( X^{\alpha}(u;A,B) \right)
\mbox{exp} \left( X^{\beta}(v;C,D) \right) :.  \label{expt}
\end{eqnarray}

\noindent It should be noticed that equations (\ref{contract})
and (\ref{expt}) hold
unchanged if we change everywhere $X^{\alpha}(u;A,B) \rightarrow
X^{\alpha}_+(u;B)$ and $X^{\beta}(v;C,D) \rightarrow X^{\beta}_-(v;C)$.

For later use let us introduce some more definitions. For $X=b^{ij},~c^{ij}$,
define

\begin{eqnarray*}
\hat{X}_\pm(u) = \mp \left( X_\pm(u;-\frac{1}{2})
- X_\pm(u; \frac{1}{2}) \right).
\end{eqnarray*}

\noindent For the bosonic fields $a^{i}(u;A,B)$, define

\begin{eqnarray}
& & \hat{a}^i_+(u) = a^i_+(u;0) - a^i_+(u; k+g), \label{ap}\\
& & \hat{a}^i_-(u) =
\frac{1}{k+g} \sum_{j,l=1}^{N-1} (B^{-1})^{jl}
\left(a^j_-(u; B_{ij}) - a^j_-(u; -B_{ij}) \right). \label{an}
\end{eqnarray}

\noindent It is easy to obtain the following operator product
expansion (OPE) relations,

\begin{eqnarray}
& & \mbox{exp} \left( \hat{a}^{i}_+(u) \right)
\mbox{exp} \left( \hat{a}^{j}_-(v) \right) \nonumber\\
& &~~~~= \frac{(u-v-B_{ij} \hbar) (u-v+(k+g+B_{ij}) \hbar)}
{(u-v+B_{ij} \hbar) (u-v+(k+g-B_{ij}) \hbar)}
\mbox{exp} \left( \hat{a}^{j}_-(v) \right)
\mbox{exp} \left( \hat{a}^{i}_+(u) \right), \label{aapn}\\
& & \mbox{exp} \left( \hat{b}^{ij}_+(u) \right)
\mbox{exp} \left( \hat{b}^{i'j'}_-(v) \right) \nonumber\\
& &~~~~= \left( \frac{(u-v)^2}
{(u-v - \hbar) (u-v+\hbar)} \right)^{\delta^{ii'} \delta^{jj'}}
\mbox{exp} \left( \hat{b}^{i'j'}_-(v) \right)
\mbox{exp} \left( \hat{b}^{ij}_+(u) \right), \nonumber\\
& & \mbox{exp} \left( \hat{c}^{ij}_+(u) \right)
\mbox{exp} \left( \hat{c}^{i'j'}_-(v) \right) \nonumber \\
& &~~~~= \left( \frac{(u-v - \hbar) (u-v+\hbar)}{(u-v)^2}
\right)^{\delta^{ii'} \delta^{jj'}}
\mbox{exp} \left( \hat{c}^{i'j'}_-(v) \right)
\mbox{exp} \left( \hat{c}^{ij}_+(u) \right).   \nonumber
\end{eqnarray}

\noindent Moreover, we have the following relations,

\begin{eqnarray*}
& & \mbox{exp} \left( \hat{b}^{ij}_+(u) \right)
: \mbox{exp} \left( \hat{b}^{i'j'}(v) \right) : \\
& &~~~~= \left( \frac{u-v-\frac{1}{2} \hbar}
{u-v +\frac{1}{2} \hbar} \right)^{\delta^{ii'} \delta^{jj'}}
: \mbox{exp} \left( \hat{b}^{i'j'}(v) \right):
\mbox{exp} \left( \hat{b}^{ij}_+(u) \right), \\
& & \mbox{exp} \left( \hat{c}^{ij}_+(u) \right)
: \mbox{exp} \left( \hat{c}^{i'j'}(v) \right) : \\
& &~~~~= \left( \frac{u-v+\frac{1}{2} \hbar}
{u-v -\frac{1}{2} \hbar} \right)^{\delta^{ii'} \delta^{jj'}}
: \mbox{exp} \left( \hat{c}^{i'j'}(v) \right):
\mbox{exp} \left( \hat{c}^{ij}_+(u) \right). \\
\end{eqnarray*}

To specify the correspondence of our notations and that of Ref. \cite{sln}
for $q$-bosons, we give the second observation

\begin{obs} \label{ob2}
The expressions $\hat{a}^{i}_\pm(u)$, $\hat{b}^{ij}_\pm(u)$
and $\hat{c}^{ij}_\pm(u)$ correspond to the fields
$a^i_{\pm}(q^{\pm \frac{k+g}{2}}z)$, $b^{ij}_\pm(z)$ and
$c^{ij}_\pm(z)$ of Ref. \cite{sln} respectively.
\end{obs}

\begin{rem} \label{rem1}
Notice that in Ref. \cite{sln}, the explicit expressions for
$a^i_{\pm}(q^{\pm \frac{k+g}{2}}z)$ are symmetric with respect to
$+ \leftrightarrow -$, but this is not the case for $\hat{a}^{i}_\pm(u)$.
The partial reason for this difference is that, for $q$-affine
algebras, the Drinfeld currents $\psi^i_\pm(z)$ are defined in an
symmetric way in $H^i_n$, whilst for Yangian doubles the currents
$H^{\pm}_i$ are defined asymmetrically.
\end{rem}

In the next subsection, we shall
see that, despite the difference stated in Remark \ref{rem1},
the above observations are
rather useful to guess the bosonic expressions for the Drinfeld
currents of $DY_\hbar(sl_N)$.

\subsection{Free boson representation of $DY_\hbar(sl_N)$ with level $k$}

Let us define

\begin{eqnarray}
& & H^{\pm}_i(u) = :
\mbox{exp} \left\{
\sum_{l=1}^i \hat{b}^{l,i+1}_{\pm}
( u \pm \frac{1}{2} (\frac{k}{2} + l -1) \hbar )
- \sum_{l=1}^{i-1} \hat{b}^{l,i}_{\pm}
( u \pm \frac{1}{2} (\frac{k}{2} + l) \hbar ) \right. \nonumber\\
& & ~~~~ + \hat{a}^i_\pm (u \mp \frac{1}{4} k \hbar ) \nonumber \\
& &~~~~ + \left. \sum_{l=i+1}^N \hat{b}^{il}_{\pm}
( u \pm \frac{1}{2} (\frac{k}{2} + l) \hbar )
- \sum_{l=i+2}^N \hat{b}^{i+1,l}_{\pm}
( u \pm \frac{1}{2} (\frac{k}{2} + l -1) \hbar ) \right\}:~, \label{hpn}
\end{eqnarray}
\begin{eqnarray}
& & E^+_i(u) = - \frac{1}{\hbar}
\sum_{m=1}^i :\mbox{exp} \left\{(b+c)^{mi}
( u + \frac{1}{2} (m-1) \hbar ) \right\} \nonumber\\
& & ~~~~\times \left[
\mbox{exp} \left( \hat{b}^{m,i+1}_+
( u + \frac{1}{2} (m-1) \hbar )
-(b+c)^{m,i+1} ( u + \frac{1}{2} m \hbar ) \right) \right. \nonumber \\
& &~~~~~~~~ - \left. \mbox{exp} \left( \hat{b}^{m,i+1}_-
( u + \frac{1}{2} (m-1) \hbar )
-(b+c)^{m,i+1}
( u + \frac{1}{2} (m-2) \hbar ) \right) \right] \nonumber\\
& &~~~~\times
\mbox{exp} \left\{ \sum_{l=1}^{m-1} \left[
\hat{b}^{l,i+1}_+ ( u + \frac{1}{2} (l-1) \hbar )
-\hat{b}^{li}_+ ( u + \frac{1}{2} l \hbar ) \right] \right\} :~, \label{ep}
\end{eqnarray}
\begin{eqnarray}
& & E^-_i(u) = - \frac{1}{\hbar} \left\{
\sum_{m=1}^{i-1} :\mbox{exp} \left( (b+c)^{m,i+1}
( u - \frac{1}{2} (k+m) \hbar ) \right) \right. \nonumber\\
& &~~~~\times \left[
\mbox{exp} \left( -\hat{b}^{mi}_- ( u - \frac{1}{2} (k+m) \hbar )
-(b+c)^{mi} ( u - \frac{1}{2} (k+m-1) \hbar ) \right) \right.  \nonumber\\
& &~~~~~~~~ - \left.
\mbox{exp} \left( -\hat{b}^{mi}_+ ( u - \frac{1}{2} (k+m) \hbar )
-(b+c)^{mi} ( u - \frac{1}{2} (k+m+1) \hbar ) \right) \right] \nonumber \\
& &~~~~\times  \mbox{exp} \left(
\sum_{l=m+1}^i \hat{b}^{l,i+1}_{-}
( u - \frac{1}{2} (k+l-1) \hbar )
- \sum_{l=m+1}^{i-1} \hat{b}^{li}_{-}
( u- \frac{1}{2} (k+l) \hbar) \right. \nonumber\\
& &~~~~~~~~ + \left.\hat{a}^i_- (u) + \sum_{l=i+1}^N \hat{b}^{il}_{-}
( u - \frac{1}{2} (k+l) \hbar)
- \sum_{l=i+2}^N \hat{b}^{i+1,l}_{\pm}
( u- \frac{1}{2} (k + l -1) \hbar ) \right): \nonumber\\
& & ~~~~+ :\mbox{exp} \left( (b+c)^{i,i+1}
( u - \frac{1}{2} (k+i) \hbar ) \right)  \nonumber\\
& & ~~~~~~~~\times \mbox{exp} \left( \hat{a}^i_- (u)
+ \sum_{l=i+1}^N \hat{b}^{il}_{-}
( u - \frac{1}{2} (k+l) \hbar )
- \sum_{l=i+2}^{N} \hat{b}^{i+1,l}_{-}
( u- \frac{1}{2} (k+l-1) \hbar) \right): \nonumber\\
& & ~~~~- :\mbox{exp} \left( (b+c)^{i,i+1}
( u + \frac{1}{2} (k+i) \hbar ) \right)  \nonumber\\
& & ~~~~~~~~\times \mbox{exp} \left( \hat{a}^i_+ (u)
+ \sum_{l=i+1}^N \hat{b}^{il}_{+}
( u + \frac{1}{2} (k+l) \hbar )
- \sum_{l=i+2}^{N} \hat{b}^{i+1,l}_{-}
( u+ \frac{1}{2} (k+l-1) \hbar) \right): \nonumber\\
& &~~~~- \sum_{m=i+2}^N :\mbox{exp} \left( (b+c)^{im}
( u + \frac{1}{2} (k+m-1) \hbar ) \right) \nonumber\\
& & ~~~~~~~~\times \left[
\mbox{exp} \left( \hat{b}^{i+1,m}_+ ( u + \frac{1}{2} (k+m-1) \hbar )
-(b+c)^{i+1,m} ( u + \frac{1}{2} (k+m) \hbar ) \right) \right. \nonumber\\
& &~~~~~~~~~~~~ - \left. \mbox{exp} \left( \hat{b}^{i+1,m}_-
( u + \frac{1}{2} (k+m-1) \hbar )
-(b+c)^{i+1,m} ( u + \frac{1}{2} (k+m-2) \hbar ) \right) \right] \nonumber\\
& &~~~~~~~~\times \left.
\mbox{exp} \left( \sum_{l=m}^{N} \left[
\hat{b}^{il}_+ ( u + \frac{1}{2} (k+l) \hbar )
-\hat{b}^{i+1,l}_+ ( u + \frac{1}{2} (k+l-1) \hbar ) \right]
\right) : \right\}~. \label{en}
\end{eqnarray}

\noindent The following proposition is the main result of this paper:

\begin{prop}   \label{prop1}
The fields $H^{\pm}_i(u)$, $E^{\pm}_i(u)$ defined in equations
(\ref{hpn}), (\ref{ep}) and (\ref{en}) are
well-defined on the Fock space ${\cal F}(l_a,l_b,l_c)$ and satisfy
equations (\ref{y1}-\ref{y4}) with $c=k$ and

\begin{eqnarray*}
& & E^{\pm}_i(u) E^{\pm}_j(v) \simeq
E^{\pm}_j(v) E^{\pm}_i(u) \sim reg. ~~\mbox{for}~B_{ij}=0,\\
& & (u-v \mp B_{ij} \hbar) E^{\pm}_i(u) E^{\pm}_j(v) \simeq
(u-v \pm B_{ij} \hbar) E^{\pm}_j(v) E^{\pm}_i(u)
\sim reg. ~~\mbox{for}~B_{ij} \neq 0,\\
& & E^{+}_i(u) E^{-}_j(v) - E^{-}_j(v) E^{+}_i(u) \\
& &~~~~ \sim reg. + \frac{1}{\hbar} \left(
\delta( u_- - v_+) H^+_i(v_+) - \delta( u_+ - v_-) H^-_i(v_-) \right),
\end{eqnarray*}

\noindent where $reg.$ means some regular expressions and $\simeq$
and $\sim$ imply ``equals up to'' such expressions.
\end{prop}

\noindent {\it Proof}: The proposition follow
by straightforward but tedious calculations. Actually the calculations
are step by step analogous to that of Ref. \cite{sln} for $q$-affine case.
So we omit all such calculations and only refer to \cite{sln} and remind the
readers of our correspondence rules (Observations 1 and 2).

\begin{rem} \label{rem2}
In proving Proposition \ref{prop1}, only the OPE relation (\ref{aapn}) for
$\hat{a}^i_\pm$ is used and the exact expressions for the fields
$\hat{a}^i_\pm$ are not important.
Actually, there are infinite many choices for
$\hat{a}^i_\pm$ which satisfy the relation (\ref{aapn}).
For example, the following
is another example which differs from the original definitions
(\ref{ap}-\ref{an}),

\begin{eqnarray}
& & \hat{a}^i_+(u) = \frac{1}{k+g} \sum_{j,l=1}^{N-1} (B^{-1})^{lj}
a^j_+(u;\frac{k+g}{2} - B_{ij})
- a^j_+(u;\frac{k+g}{2} + B_{ij}), \label{ap2}\\
& & \hat{a}^i_-(u) = a^i_-(u;-\frac{k+g}{2})
- a^i_+(u;\frac{k+g}{2}). \label{an2}
\end{eqnarray}

\noindent However, no matter which choice we use, we cannot make the
definition of $\hat{a}^i_+(u)$ and $\hat{a}^i_-(u)$ symmetric, i.e.
no violation of Remark \ref{rem1} could occur.
\end{rem}

\begin{rem} \label{rem3}
While $N=2$, equations (\ref{hpn}-\ref{en}) become

\begin{eqnarray}
& & H^{+}(u) = : \mbox{exp}\left( \hat{b}_+(u+\frac{1}{4}k\hbar) +
\hat{b}_+(u+\frac{1}{2}(\frac{k}{2}+2) \hbar)
+ \hat{a}_+(u-\frac{1}{4}k\hbar) \right):~, \label{hppp}\\
& & H^{-}(u) = : \mbox{exp}\left( \hat{b}_-(u-\frac{1}{4}k\hbar) +
\hat{b}_-(u-\frac{1}{2}(\frac{k}{2}+2) \hbar)
+ \hat{a}_-(u+\frac{1}{4}k\hbar) \right):~,\label{hnnn}\\
& & E^+(u)= -\frac{1}{\hbar} :  \left[
\mbox{exp}\left( \hat{b}_+(u) -(b+c)(u+\frac{\hbar}{2}) \right)
- \mbox{exp}\left( \hat{b}_-(u) -(b+c)(u-\frac{\hbar}{2})
\right) \right]:~, \label{eppp}\\
& & E^-(u)= \frac{1}{\hbar} :  \left[
\mbox{exp}\left( (b+c)(u+\frac{1}{2}(k+1)\hbar ) \right)
\mbox{exp}\left( \hat{a}_+(u) + \hat{b}_+(u+ \frac{1}{2}(k+2)\hbar)
\right) \right. \nonumber \\
& &~~~~\left.
- \mbox{exp}\left( (b+c)(u-\frac{1}{2}(k+1)\hbar ) \right)
\mbox{exp}\left( \hat{a}_-(u) + \hat{b}_-(u- \frac{1}{2}(k+2)\hbar)
\right) \right]:  \label{ennn}
\end{eqnarray}

\noindent which is different from the result of Ref. \cite{konno} for
$DY_\hbar(sl_2)_k$. The reason for this difference is that, first, the
Yangian double $DY_\hbar(sl_2)_k$ of Ref. \cite{konno}
is realized in an asymmetric
way which differs from the symmetric one which we are using by a
shift of parameter \cite{KL};
Second, as we have remarked in Remarks \ref{rem1}
and \ref{rem2}, for the same realization of Yangian double, there still
exist infinite many choices for the bosonization formulas. Therefore
the difference of our result (\ref{hppp}-\ref{ennn}) from that of
Ref. \cite{konno} is reasonable.
\end{rem}

\section{Problems in obtaining Vertex operators and screening currents}

After successfully obtained the bosonization formulas for the Yangian
double $DY_\hbar(sl_N)$ using the correspondence rules (Observations
\ref{ob1} and \ref{ob2}), one naturally expects that the Vertex operators
and screening currents for  $DY_\hbar(sl_N)$ could also be obtained
in the same way. In this section we briefly give why this is difficult.

Following Observations \ref{ob1} and \ref{ob2} and Ref.\cite{sln}, we expects
that the screening currents for $DY_\hbar(sl_N)$ might be written in the
following form,

\begin{eqnarray*}
S^i(u) = :\mbox{exp}\left( X^i [a](u) \right): \tilde{S}^i(u),
\end{eqnarray*}

\noindent where $\tilde{S}^i(u)$ is nothing but $E^+_{N-i}(u)$ with the
replacement $\hat{b}^{ij}_\pm \rightarrow - \hat{b}^{N+1-j,N+1-i}_\mp$,
$(b+c)^{ij} \rightarrow (b+c)^{N+1-j,N+1-i}$, and $X^i [a](u)$
is some field depending only on $a^i$ but not on $b^{ij}$ and $c^{ij}$.
In the $q$-affine case, $X^i [a](u)$ is just the field
$-\left( \frac{1}{k+g} a^i \right)(z; \frac{k+g}{2})$ \cite{sln}.
At present, we expect that $:\mbox{exp}\left( X^i [a](u) \right):$
have the following OPE relations with
$\mbox{exp}\left( \hat{a}^i_+(u) \right)$
and $\mbox{exp}\left( \hat{a}^i_-(u) \right)$ (see equations (C.17), (C.18)
of Ref. \cite{sln})

\begin{eqnarray}
& & \mbox{exp} ( \hat{a}^i_+(u) )
:\mbox{exp} ( X^j [a](v) ): \nonumber \\
& &~~~~~~= \frac{u-v+(\frac{k+g}{2}-B_{ij})\hbar}
{u-v+(\frac{k+g}{2}+B_{ij})\hbar}
:\mbox{exp} ( X^j [a](v) ):
\mbox{exp} ( \hat{a}^i_+(u) ),  \label {1a}\\
& & :\mbox{exp} ( X^i [a](u) ):
\mbox{exp} ( \hat{a}^j_-(v) ) \nonumber \\
& &~~~~~~= \frac{u-v+(\frac{k+g}{2}-B_{ij})\hbar}
{u-v+(\frac{k+g}{2}+B_{ij})\hbar}
\mbox{exp} ( \hat{a}^j_-(v) )
:\mbox{exp} ( X^i [a](u) ):~.  \label{2a}
\end{eqnarray}

\noindent Notice that in the rational factors of equations
(\ref{1a}) and (\ref{2a}),
the $\frac{k+g}{2}$ appear with the same sign in both the numerators and
the denominators. This fact makes it difficult to obtain an explicit
expression for $X^i [a](u)$. For examples, if we adopt the definitions
(\ref{ap}) and (\ref{an}) of $\hat{a}^i_\pm(u)$,
then the ``positive frequency part''
of $X^i [a](u)$ can be easily seen to be equal to
$\hat{a}^i(u+\frac{k+g}{2})$, but the ``negative frequency part
could not be written in a simple form (and it is not known whether it is
possible to write down such an expression). If we adopt the definitions
(\ref{ap2}) and (\ref{an2}) instead of (\ref{ap}) and (\ref{an}),
then the negative frequency part of $X^i [a](u)$
can be obtained easily but the positive frequency part is unknown.
That is why we could not obtain a simple analogy of screening currents
for Yangian double and $q$-affine algebras. Due to similar reasons
the vertex operators for $DY_\hbar(sl_N)$ is also not obtained from that of
$U_q(\widehat{sl_N})$.

\section{Discussions}

In this paper we established the free boson representation of the Yangian
double $DY_\hbar(sl_N)$ with arbitrary level $k$. Our construction
is based on the crucial correspondence Observations \ref{ob1} and \ref{ob2}.
Such representations of the Yangian double $DY_\hbar(sl_N)$ are
expected to be useful in calculating the correlation functions
of various quantum integrable systems in (1+1)-spacetime dimensions,
e.g. the spin Calogero-Sutherland model \cite{spinCS}, quantum nonlinear
Schrodinger equation \cite{nonlinear} and some field theoretic
models such as Thirring model, Gross-Neveu model with $U(N)$
gauge symmetries etc.
Our representation of $DY_\hbar(sl_N)$ may also be used to analyze the
behavior of Yangian double at the critical level $k=-g$, a very fascinating
area of great interest of study \cite{FR}.

Besides what have been solved in this paper, the unsolved problem of
the construction of vertex operators and screening currents are also
of great interest. Especially if we know these quantities we could have
been able to calculate the cohomology of the action of our
bosonization formulas on the Fock spaces ${\cal F}(l_a,~l_b,~l_c)$.
We hope these problems could be solved in the future.

\end{document}